\documentclass[12pt,a4paper]{article}
\usepackage{amsfonts,amsmath,amssymb,amscd}
\usepackage{fancybox}
\usepackage[latin1]{inputenc}
\usepackage[usenames]{color}
\usepackage{graphicx}
\usepackage{subfig}
\usepackage{float}
\usepackage{fullpage}
\usepackage{epsfig,amssymb}

\newtheorem{thm}{Theorem}[section]

\newtheorem{prop}[thm]{Proposition}

\DeclareMathOperator{\arcsinh}{arcsinh}
\DeclareMathOperator{\arccosh}{arccosh}

\begin{document}
	
\title{\LARGE \textbf{Visual Curve Completion and Rotational Surfaces of Constant Negative Curvature}}

\author{\Large{\'Alvaro P\'ampano}\vspace{0.5cm}\\Department of Mathematics\\Faculty of Science and Technology\\University of the Basque Country\\Bilbao (Spain)\\\textit{alvaro.pampano@ehu.es}}
\date{\today}

\maketitle

\begin{abstract}
If a piece of the contour of a picture is missing to the eye vision, then the brain tends to complete it using some kind of sub-Riemannian geodesics of the unit tangent bundle of the plane, $\mathbb{R}^2\times\mathbb{S}^1$. These geodesics can be obtained by lifting extremal curves of a total curvature type energy in the plane.\\ 
We completely solve this variational problem, geometrically. Moreover, we also show a way of constructing rotational surfaces of constant negative curvature in $\mathbb{R}^3$ by evolving these extremal curves under their associated binormal flow with prescribed velocity.\\
Finally, we prove that, locally, all rotational constant negative curvature surfaces of $\mathbb{R}^3$ are foliated by extremal curves of these energies. Therefore, we conclude that there exists a one-to-one correspondence between the sub-Riemannian geodesics used by the brain for visual curve completion and these rotational surfaces of $\mathbb{R}^3$.\\

\noindent{\emph{Keywords:} Constant Gaussian Curvature Surfaces · Extremal Curves · Sub-Riemannian Geodesics · Total Curvature Type Energy · Unit Tangent Bundle.}
\end{abstract}

\section{Introduction}

These notes are a printed version of the talk given by the author at the 23rd International Summer School on Global Analysis and Applications held in Brasov in August 2018. The purpose of the talk was to present some results included in the works \cite{AGP1}, \cite{AGP2} and \cite{P}. Here, ideas and arguments are only sketched while proofs are omitted. Interested readers are going to be referred to \cite{AGP1}, \cite{AGP2}, or \cite{P}, respectively, for a complete and more general treatment.

In this paper we aim to explain a variational approach for visual curve completion purposes, which is supposed to imitate the mechanisms of our brain. For different approaches one may see \cite{DBRS}.

Neuro-biologic research over the past few decades has great\-ly clarified the functional mechanisms
of the first layer V1 of the visual cortex (\textit{primary visual cortex}). Such layer contains a variety of types of cells, in\-clu\-ding
the so-called simple cells. Researchers found that V1 constitutes of orientation selective cells at all o\-rien\-tations for all retinal positions so simple cells are sensitive to orientation specific brightness gradients (for details see \cite{BYBS} and \cite{DBRS}). Recently, this structure of the primary visual cortex has been modeled using  sub-Riemannian geometry, \cite{Petitot}. In particular, the unit tangent bundle of the plane can be used as an abstraction to study the organization and mechanisms of V1.

It is believed that, if a piece of the contour of a picture is missing to the eye vision (or maybe it is co\-ve\-red by an object), then the brain tends to ``complete" the curve by minimizing some kind
of energy, be\-ing length the simplest (but not the only) of such (for other possible options see \cite{AGP2}). In short, \emph{there is some sub-Riemannian structure on the space of visual cells and the brain considers a sub-Riemannian geodesic between the endpoints of the missing data}.

As pointed out in, for instance, \cite{AGP2}, \cite{BYBS} and \cite{DBRS}, it turns out that these sub-Riemannian geodesics can be projected down to critical curves of 
\begin{equation}
\mathbf{\Theta}(\gamma)=\int_\gamma \sqrt{\kappa^2+1}\, ds\, \label{1}
\end{equation}
in the Euclidean plane, $\mathbb{R}^2$.

In order to solve this variational problem for fixed boundary conditions we deal with a highly nonlinear system that must be solved. Therefore, due to this difficulty, it seems reasonable to develop a numerical approach. In \S 2, we describe this model for visual curve completion and show how the XEL-platform introduced in \cite{AGP2} for this particular case can be used to numerically complete damaged or covered images.

However, due to the hypercolumnar organization of visual cells, it is more accurate to consider critical curves of the following curvature type energy
\begin{equation}
\mathbf{\Theta}_a(\gamma)=\int_\gamma \sqrt{\kappa^2+a^2}\, ds\, \label{3}
\end{equation}
acting on planar non-geodesic curves of $\mathbb{R}^2$. Therefore, in \S 3 we completely solve this variational problem, geometrically. That is, we explicitly obtain the curvature of the critical curves and, as it is well-known, these curvatures completely describe the curves, up to rigid motions.

Notice that, although the sub-Riemannian model of \S 2 for visual curve completion only considers curves, we are going to show the close relation between critical curves of $\mathbf{\Theta}_a$ and a well-known family of rotational surfaces of $\mathbb{R}^3$ with a nice geometric property. Indeed, in \S 4 we variationally characterize the profile curve of rotational surfaces with constant negative Gaussian curvature as critical curves of $\mathbf{\Theta}_a$, \eqref{3}. Moreover, we also describe the binormal evolution of these extremal curves, giving rise to a natural construction of constant negative curvature rotational surfaces.

In this way, we prove the existence of a one-to-one correspondence between sub-Riemannian geodesics of the unit tangent bundle of the plane, $\mathbb{R}^2\times\mathbb{S}^1$, used by the brain for visual curve completion and rotational surfaces with constant negative Gaussian curvature in the Euclidean 3-space, $\mathbb{R}^3$. This correspondence is quite natural in the sense that it comes intrinsically from the properties of the projections into the plane of the corresponding sub-Riemannian geodesics. Therefore, it suggests that the model based on the use of a sub-Riemannian structure in the unit tangent bundle of the plane in order to study the functional mechanisms of the cortex V1 may give not only the completed curve, but also, some extra information as we explain in \S 5.

\section{Visual Curve Completion}

The main objective of this section consists on introducing a problem concerning visual curve completion and showing how it is related with a variational problem over curves where the corresponding energy depends on the curvature. 

\begin{figure}[H]
\begin{centering}
{\includegraphics[width=5.5cm,height=5.5cm]{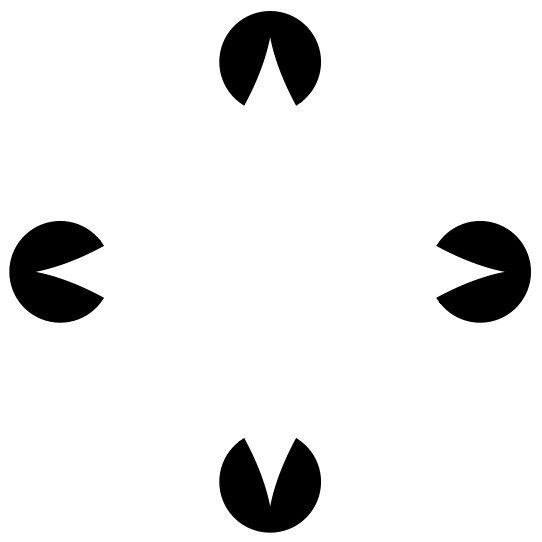}}\hspace{2cm}{\includegraphics[width=5.5cm,height=5.5cm]{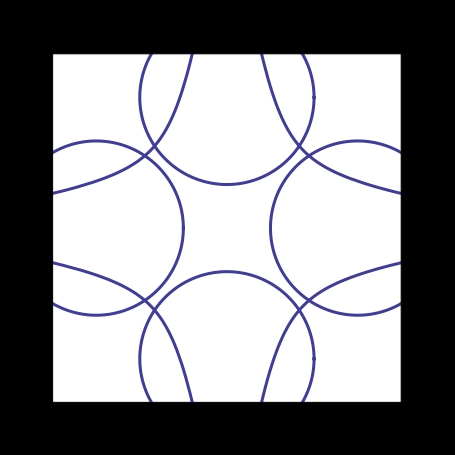}}
\caption{Two kind of covered or damaged images.}\label{s}
\end{centering}
\end{figure}

Assume that we want to recover an image that is covered or damaged, as it is represented in Figure \ref{s}. Here, we can see two different images that are not complete for different reasons. This is a case that often appears in real world, due to the many natural or artificial obstacles that may appear when, for instance, taking a photo. Moreover, some information of real photos may be lost, which gives rise, once more, to damaged images. Therefore, our purpose consists on developing a way of completing these images in a natural way, that is, trying to imitate the intuitive answer our brain plots.

Then, according to the model of \cite{P}, in the space $\mathbb{R} ^2 \times  \mathbb{S} ^1$ (which models the layer V1) each point $(x, y, \theta)$ represents a column of cells associated to a point of \emph{retinal data} $(x, y) \in \mathbb{R} ^2$, all of which are adjusted to the \emph{orientation} given by the angle $\theta\in \mathbb{S}^1$. In other words, the vector $(\cos\theta, \sin\theta)$ is the direction of maximal rate of change of brightness at point $(x, y)$
of the picture seen by the eye. Such vector can be seen as the normal to the boundary of the picture. Thus, when the cortex cells are stimulated by an image, the border of the image gives
a curve inside the $3$-space $\mathbb{R} ^2 \times  \mathbb{S} ^1$, but such curves are restricted to be tangent to the distribution spanned by the vector fields
\begin{equation}
X_1 = \cos\theta\, \frac{\partial}{\partial x} + \sin\theta\,\frac{\partial}{ \partial y} \, , \quad\quad X_2 = \frac{\partial}{ \partial \theta}\, \label{2.25}.
\end{equation}

Therefore, the existence of this distribution implies that we need to define an adequate sub-Riemannian structure on $\mathbb{R}^2\times\mathbb{S}^1$. First, we take the distribution $\mathcal{D} = ker(\sin\theta  \, dx - \cos\theta\, dy)$, where $x$ and $y$ are the coordinates on $\mathbb{R}^2$ and $\theta$ is the coordinate on $\mathbb{S}^1$. This distribution is spanned by the vector fields described in \eqref{2.25}. Then, we consider on $\mathcal{D}$ the inner product $\langle\cdot ,\cdot \rangle$ for which the  two vectors \eqref{2.25} are everywhere orthonormal. Finally, it is believed that our brain minimizes some kind of energy in this ambient space to complete the images. Along this paper, we are going to consider the sub-Riemannian geodesic between the endpoints of the missing data as the completed curve (for other possible choices see \cite{AM}, \cite{AGP2} and \cite{BYBS2}).

\begin{figure}[H] 
\begin{centering}
{\includegraphics[width=5.5cm,height=5.5cm]{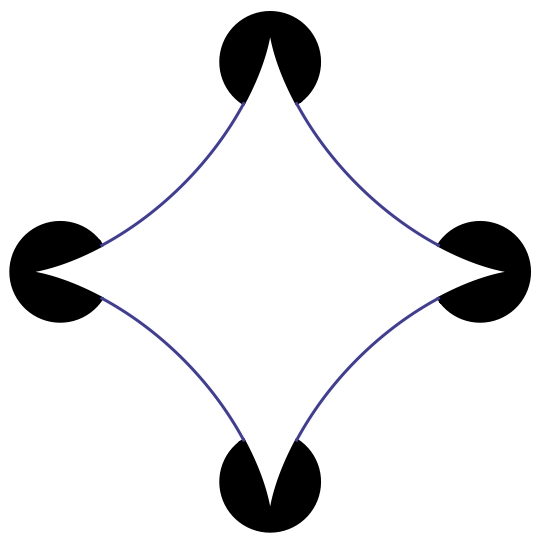}}\hspace{2cm}{\includegraphics[width=5.5cm,height=5.5cm]{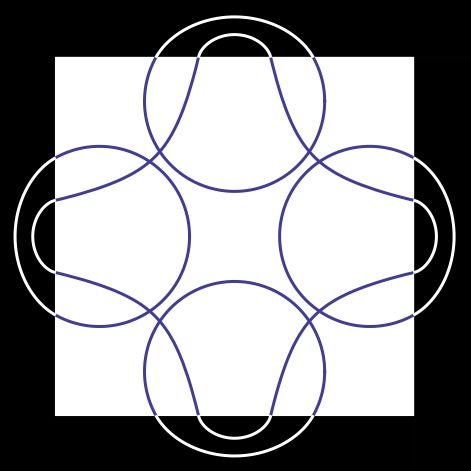}}
\caption{Completions of Figure \ref{s}.}\label{con}
\end{centering}
\end{figure}

In Figure \ref{con}, we have the completed images corresponding with Figure \ref{s}, where we have used a direct numerical approach based on the gradient descent method implemented in the XEL-platform. For more details about the XEL-platform and the algorithm based on the gradient descent method, see \cite{AGP2}.

Now, notice that every $\mathcal{D}$-curve in $\mathbb{R}^2\times\mathbb{S}^1$ (that is, a smooth immersed curve which is always tangent to $\mathcal{D}$) $\gamma(t) = (x(t), y(t),$ $\theta(t))$ with $\gamma^*(\cos\theta  \, dx +\sin\theta\, dy)\neq 0$ is the lift of a regular curve $\alpha(t) = (x(t), y(t))$ in the plane whose tangent vector $\alpha ' (t)$ forms the angle $\theta(t)$ with the $x$-axis, that is,
\begin{equation}
\alpha ' (t) = v(t)\, \cos\theta\, \frac{\partial}{\partial x} + v(t)\, \sin\theta\,\frac{\partial}{ \partial y} \, ,\nonumber
\end{equation}
where $v(t)$ is the speed of $\alpha(t)$. Conversely, every regular curve $\alpha(t)$ in the plane may be lifted to a $\mathcal{D}$-curve $\gamma(t) = (x(t), y(t), \theta(t))$ by setting
$\theta(t)$ equal to the angle between $\alpha '(t)$ and the $x$-axis. Now, the tangent vector $\gamma'(t)$  of the $\mathcal{D}$-curve $\gamma(t) $ has squared length
\begin{eqnarray}
\langle \gamma ' (t), \gamma ' (t)\rangle  =  v^2(t)+\left(\theta'\right)^2(t)  =  v^2(t)\left(1+\left(\frac{\theta'(t)}{v(t)}\right)^2\right) = v^2(t)(1+\kappa^2(t))  \,  ,\nonumber
\end{eqnarray}
where $\kappa(t)$ is the curvature of $\alpha$. Thus, the $\mathcal{D}$-curves with $\gamma^*(\cos\theta  \, dx +\sin\theta\, dy) \neq 0$ that cover the distance between two points
$(x_0, y_0, \theta_0)$ and $(x_1, y_1, \theta_1)$ of $\mathbb{R}^2\times\mathbb{S}^1$ are the lifts of curves $\alpha$ in the plane joining $(x_0, y_0)$  to $(x_1, y_1)$
with initial angle $\theta_0$ and final angle $\theta_1$ that minimize the functional $\mathbf{\Theta}$, \eqref{1}, among all such curves in the plane. In other words, \emph{geodesics in V1 are obtained by lifting to $\mathbb{R}^2 \times  \mathbb{S}^1$ minimizers of \eqref{1} in $\mathbb{R}^2$}. 

Finally, as indicated in \cite{BYBS}, the hypercolumnar organization of the visual cortex suggests that the cost of moving one orientation unit is not necessarily the same as to moving spatial units, then the curve completion problem should consider the functional $\mathbf{\Theta}_a$, \eqref{3}, for any real constant $a$, ac\-ting on planar curves instead. 

\section{Total Curvature Type Energies in the Plane}

Let us consider along this section the energy functional $\mathbf{\Theta}_a \colon
\Omega_{pq} \rightarrow \mathbb R $ defined in \eqref{3} where, $a\in \mathbb{R}$, $s$ is  the arc-length parameter and $ \kappa (s)$ is the curvature of $\gamma (s)$, acting on the space of immersed curves $\gamma (t)$ in  $\mathbb{R}^2 $ with fixed endpoints $p,q\in \mathbb{R}^2$, $\Omega_{pq}$. 

Observe that geo\-de\-sics are minima of $\mathbf{\Theta}_a $, \eqref{3}. Now, since geodesics are always critical  for $\mathbf{\Theta}_a$, \eqref{3} (for suitable boundary conditions) in the following we may assume, in addition, that $\gamma\in \Omega_{pq}$ is a non-geodesic curve. Then, using standard arguments that involve integration by parts, we obtain that the Eu\-ler-Lagrange operator is given by
\begin{eqnarray}
\mathcal{E}(\gamma)& =& \frac{1}{\left(\kappa^{2}+a^2\right)^{\frac{1}{2}}}\,T^{(3)}+ 2\frac{d}{ds}\left(\frac{1}{\left(\kappa^{2}+a^2\right)^{\frac{1}{2}}}\right)T''\nonumber\\&&+ \left(\frac{d^2}{ds^2}\left(\frac{1}{\left(\kappa^{2}+a^2\right)^{\frac{1}{2}}}\right)+ \frac{\kappa^{2}-a^2}{\left(\kappa^{2}+a^2\right)^{\frac{1}{2}}} \right)T'+  \frac{d}{ds}\left(\frac{\kappa^{2}-a^2}{\left(\kappa^{2}+a^2\right)^{\frac{1}{2}}}\right) T\,  ,\nonumber
\end{eqnarray}
$T$ being the unit tangent vector field to $\gamma$. Now, if $\gamma$ is an extremal of  $\mathbf{\Theta}_a$, \eqref{3}, then $\mathcal{E}(\gamma)=0$ and, using above equation we get the following \emph{Euler-Lagrange equation}
\begin{equation}
\frac{d^2}{ds^2}\left(\frac{\kappa}{\left(\kappa^{2}+a^2\right)^{\frac{1}{2}}}\right) - a^2\, \frac{\kappa}{ \left(\kappa^{2}+a^2\right)^{\frac{1}{2}} }= 0 \,  .\label{EL}
\end{equation}
Case $a=0$ corresponds to the\textit{ total curvature functional}. This is the reason why along this paper we call \emph{total curvature type energy} to $\mathbf{\Theta}_a$, \eqref{3}. In the proper total curvature functional case, if $\gamma\in \Omega_{pq}$ is critical then,
\begin{equation}
\int_\gamma \kappa = \theta (1)-\theta(0)+ 2 \pi \,m \, , \label{2.34}
\end{equation}
where $\theta (i)\in[0,2\pi), i=1,2,$ denotes the angle that $v_i, i=1,2,$ makes with the line determined by $p$ and $q$, and $m\in \mathbb{Z}$ is an integer representing the number (taking orientation into account)  of loops that the trace of $\gamma$ has  between $p$ and $q$. Since \eqref{2.34} is constant within any regular homotopy class of $\Omega_{pq}$ having the same tangent vectors at $p$ and $q$, we see that the corresponding variational problem is trivial, that is, any $\gamma \in\Omega_{pq}$ with adequate tangent vectors at $p$ and $q$ is critical for $\mathbf{\Theta}_a$, \eqref{3}, when $a=0$.

So we assume from now on that $a\neq0 $. Notice that above equation \eqref{EL} imply that there are no critical circles for $\mathbf{\Theta}_a$, \eqref{3}. On the other hand, if $\gamma\in\Omega_{pq}$ is critical for $\mathbf{\Theta}_a$ under arbitrary boundary conditions, then it satisfies $\mathcal{E}(\gamma)=0$ and, therefore, we conclude that the first integrals of the Euler-Lagrange equation, \eqref{EL}, is
\begin{equation}
\left(\frac{d\kappa}{ds}\right)^2=\frac{\left(\kappa^2+a^2\right)^2}{a^4} \left(\,d\,(\kappa^2+a^2) - a^4\right)\,  , \label{2.35}\\
\end{equation}
where $d>0$ is a constant of integration.

We recall now that a vector field $W$ along $\gamma$, which infinitesimally preserves unit speed parametrization is said to be a \emph{Killing vector field along $\gamma$} (in the sense of \cite{Langer-Singer-4}) if $\gamma$ evolves in the direction of $W$ without changing shape, only position. In other words, the following equations must hold
\begin{equation}
W(v)(s,0)=W(\kappa)(s,0)=0\, , \label{W}
\end{equation}
($v=\lvert \gamma'\rvert$ being the speed of $\gamma$) for any variation of $\gamma$ having $W$ as variation field.

It turns out that critical curves of $\mathbf{\Theta}_a$, \eqref{3}, in $\mathbb{R}^2$ have a naturally associated Killing vector field defined along them. Let us define the following vector field along $\gamma$
\begin{equation}
\mathcal{J}=\frac{d}{ds}\left(\frac{\kappa}{\left(\kappa^2+a^2\right)^\frac{1}{2}}\right)N-\,\frac{a^2}{\left(\kappa^2+a^2\right)^\frac{1}{2}}\,\,T\, , \label{j}
\end{equation}
where $N$ represents the unit normal vector field to the curve $\gamma$.

Now, as proved for instance in \cite{Langer-Singer-4}, combining both the Euler-Lagrange equation of a planar critical curve, \eqref{EL}, and the definition of Killing vector fields along curves, \eqref{W}, we obtain

\begin{prop}\label{killing} 
Assume that $\gamma$ is a critical curve of $\mathbf{\Theta}_a$, \eqref{3}, in $\mathbb{R}^2$ and consider the vector field along $\gamma$, $\mathcal{J}$ defined in \eqref{j}. Then $\mathcal{J}$ is a Killing vector fields along $\gamma$.
\end{prop}

Moreover, by using the symmetries of $\mathbb{R}^2$ and above Killing vector field along critical curves, the coordinates of these curves $\gamma$ can also be obtained by quadratures. In fact, $\gamma$ being critical means that $\mathcal{J}$, \eqref{j}, can be extended to a Killing field defined on the whole $\mathbb{R}^2$ (see \cite{Langer-Singer-4}, for more details). We denote this extensions with the same letter, $\mathcal{J}$. Furthermore, this extension must come from a one-parameter group of motions. Indeed, using equation \eqref{2.35} we see that the length of the Killing field $\mathcal{J}$ is constant, $\lvert \mathcal{J}\rvert=d$. Therefore, $\mathcal{J}$ must be a translational vector field. That is, choosing $z$ as the direction of the motion associated to $\mathcal{J}$, introducing coordinates $r$, $z$, and taking into account the first integral of the Euler-Lagrange equation given before, \eqref{EL}, we get (one may also see \cite{Langer-Singer-2})
\begin{equation}
\mathcal{J}= \sqrt{d} \, \partial_z \,  , \label{2.37}
\end{equation}
where $d>0$ is the constant of integration appearing in \eqref{2.35}. Now, since $\gamma(s)=\left(r(s),z(s)\right)$, we have that $T(s)=r'\partial _r+ z' \partial_z$ so \eqref{j} and \eqref{2.37} imply
\begin{equation}
z'(s)=\frac{-a^2}{\left(d(\kappa^2+a^2)\right)^\frac{1}{2}} \, . \label{2.40}
\end{equation}
Moreover, assuming that $s$ represents the arc-length parameter and using above equation, \eqref{2.40}, we conclude that
\begin{equation}
r^2(s)=\frac{\kappa^2}{d(\kappa^2+a^2)}\, . \label{2.39}
\end{equation}

Therefore, from \eqref{2.40} and \eqref{2.39} we see that, once the curvature of the critical curve, $\gamma$, is known (for which we need to solve \eqref{2.35}) the coordinates of $\gamma$ can be obtained by quadratures. Notice also that \eqref{2.40} implies that $z(s)$ is monotonic, thus, there are no periodic solutions of the Euler-Lagrange equation \eqref{EL}, so that we cannot have closed critical curves.

Moreover, extremal curves in $\mathbb{R}^2$ are totally determined by their curvature, $\kappa$, which, in our case, can be obtained explicitly. In fact, from the Euler-Lagrange equation, \eqref{EL}, we have

\begin{prop}\label{curvaturas} 
Let $\gamma\subset\mathbb{R}^2$ be a critical curve of $\mathbf{\Theta}_a$, \eqref{3}. Then, either $\gamma$ is a geodesic or $\mathbf{\Theta}_a$, \eqref{3}, is acting on $\Omega_{pq}$ and the curvature of $\gamma$ is given by
\begin{equation}
\kappa^2(s)=\frac{a^2\left(d-a^2\right)f^2\left(a\,s\right)}{a^2-\left(d-a^2\right)f^2\left(a\,s\right)}\, ,\label{2.42}
\end{equation}
where $f(z)=\sinh z$, $\cosh z$ or $e^z$.
\end{prop}

From this equation, \eqref{2.42}, we see that there is at most one point where the curvature may change sign, so the curvature is always positive (or negative) except for at most one inflection point.  Moreover $\kappa' (s)$ has at most one zero (a vertex), hence, either the curvature is monotonic, or it monotonically decreases up to reaching the vertex where it starts to monotonically increase (or viceversa). This means that planar critical  curves are a family of ``short time" spirals in the plane.

This is a case relevant in i\-ma\-ge restoration as we mentioned in the introduction and a parametrization of extremal curves, using as parameter, $\theta$,  the angle which the curve makes with a fixed line, was given in \cite{BYBS}, under the assumption that the curves have no inflection points. But, as we have just noticed, these critical curves have at most one vertex, hence, the argument of \cite{BYBS} applies and an explicit parametrization of extremal curves for this variational problem can be obtained in terms of elliptic integrals of the first and second kind (see \cite{AGP2} for details). Alternatively, one may use our previous computations to get different parametrizations of extremals, at least, by quadratures. In fact, as it is very well-known, a parametrization in terms of the curvature and arc-length parameter of a planar curve is given by $\left(\int\cos\int\kappa,\int\sin\int\kappa \right)$, then, using \eqref{2.42} we can also get a parametrization of the extremals of $\mathbf{\Theta}_a$ in $\mathbb{R}^2$ in terms of the arc-length parameter after two quadratures. Observe that another parametrization can also been explicitly obtained from \eqref{2.40} and \eqref{2.39}.

\begin{figure}[H] 
\begin{centering}
{\includegraphics[width=4.5cm,height=4.5cm]{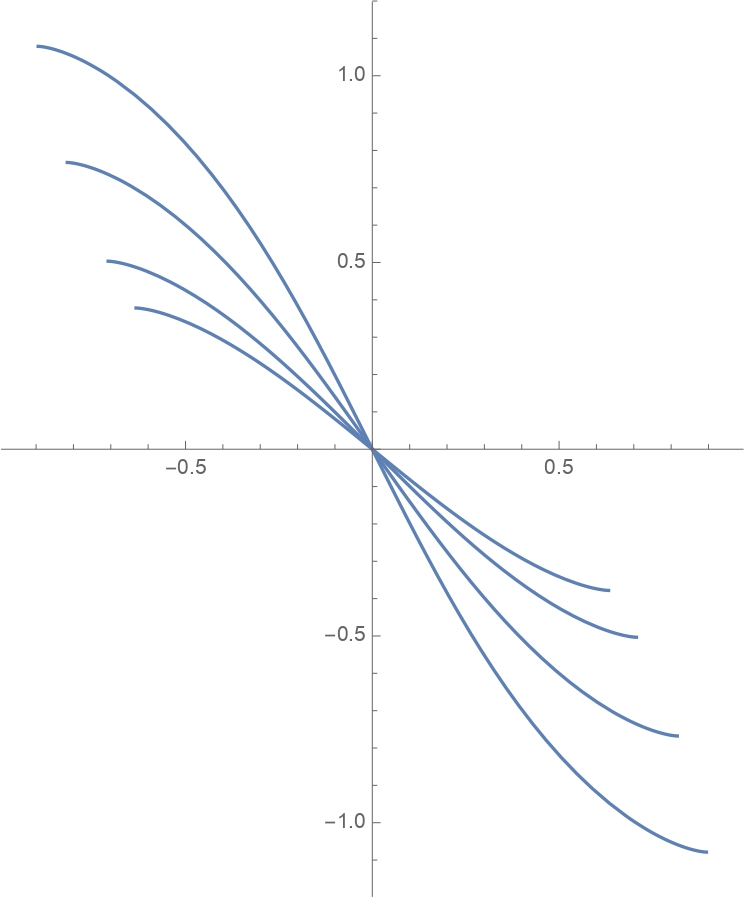}}\hspace{1cm}{\includegraphics[width=4.5cm,height=4.5cm]{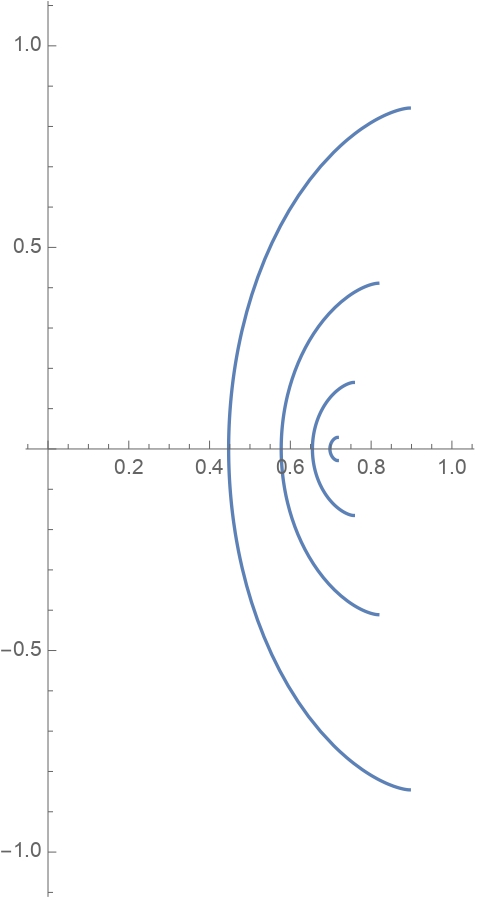}}\hspace{1cm}{\includegraphics[width=4.5cm,height=4.5cm]{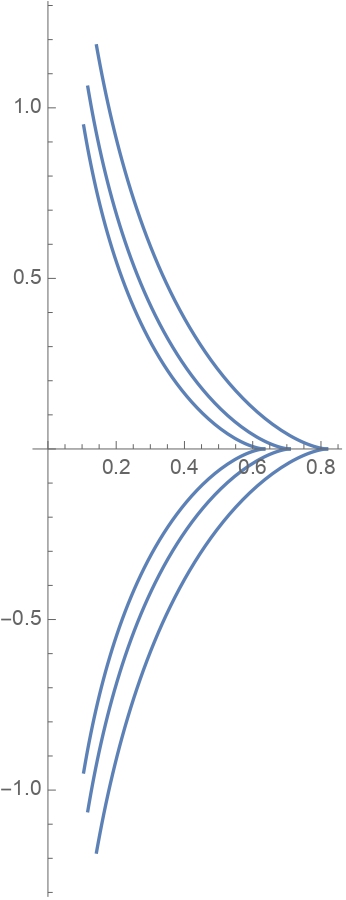}}
\caption{Critical curves of $\mathbf{\Theta}_a$ with curvatures given in \eqref{2.42} for $f(z)=\sinh z$ (Left), $f(z)=\cosh z$ (Center) and $f(z)=e^z$ (Right).}
\end{centering}
\label{figure1}
\end{figure}

However, for any possible choice of a parametrization method, a specific determination of the solution curves implies that the integration constants must be determined. This can be tried by imposing the solutions to satisfy the given boundary conditions, but this requires, at the best, solving a highly nonlinear system for which an explicit parameters expression looks unlikely. Hence, a numerical approach seems to be a reasonable strategy, as we have explained in \S 2.

\section{Rotational Surfaces with Constant Negative Curvature}

Along this section we are going to consider that $\gamma\subset\mathbb{R}^2$ is any non-geodesic critical curve of $\mathbf{\Theta}_a$, \eqref{3}, acting on $\Omega_{pq}$. Moreover, for the purpose of this section, we are going to consider that $\mathbb R^2$ is contained in the Euclidean 3-space, $\mathbb R^3$. In this setting, extremal curves have another naturally associated Killing vector field in the direction of the binormal which can be uniquely extended to a Killing vector field in the whole $\mathbb{R}^3$. Therefore, in this section, we are going to evolve these critical curves under their associated Killing vector field in the direction of the binormal. 

Let us define the vector field along $\gamma$ in the direction of the binormal by
\begin{equation}
\mathcal{I}=\frac{\kappa}{\left(\kappa^2+a^2\right)^\frac{1}{2}}\, B\,, \label{i}
\end{equation}
where $\kappa$ denotes the curvature of the critical curve $\gamma$ and $B$ is the unit binormal vector field.

Now, a similar argument as the one applied in Proposition \ref{killing} proves that equations \eqref{W} are verified for $W=\mathcal{I}$ and, therefore, $\mathcal{I}$ is a Killing vector field along $\gamma$. Then, following the arguments of \cite{Langer-Singer-4}, this vector field along $\gamma$, $\mathcal{I}$, can be uniquely extended to the whole $\mathbb R^3$. We are going to use $\xi$ to represent this extension and we are going to denote by $\{\phi_t\, ;\, t\in \mathbb R\}$ the one-parameter group of isometries associated to $\xi$, that is, the flow of $\xi$. 

Now, we locally define the \emph{binormal evolution surface} (for more details see \cite{AGP1} and references therein)
\begin{equation}
S_\gamma:=\{ x(s,t)=\phi_t\left(\gamma(s)\right)\}\, . \nonumber
\end{equation}
Then, we have the following result

\begin{prop}[\cite{AGP2}] 
Let $\gamma$ be a critical curve of $\mathbf{\Theta}_a$, \eqref{3}, acting on $\Omega_{pq}$. Then, the binormal evolution surface with initial condition $\gamma$ defined above is a rotational surface.
\end{prop}

Notice that planar extremal curves of $\mathbf{\Theta}_a$, \eqref{3}, have non-constant curvature and, as a consequence, their associated binormal evolution surface is not flat. However, the following result about their Gaussian curvature can be proved

\begin{thm}[\cite{P}] \label{thm1}
Let $\gamma$ be a planar non-geodesic critical curve of $\mathbf{\Theta}_a$ acting on $\Omega_{pq}$. Then, the binormal evolution surface with initial condition $\gamma$, $S_\gamma$, is a rotational surface with constant negative Gaussian curvature. Indeed, the Gaussian curvature of $\gamma$ is given by $K=-a^2$.
\end{thm}

Using Theorem \ref{thm1} we have a natural way of constructing rotational surfaces with fixed constant negative curvature by evolving planar critical curves of $\mathbf{\Theta}_a$, \eqref{3}, under their associated binormal flow with prescribed velocity $\lvert\mathcal{I}\rvert$, see \eqref{i}.

Moreover, the converse of Theorem \ref{thm1} is also true. Indeed, in what follows we are going to show that for any rotational surface with negative constant Gaussian curvature $K<0$, it is possible to find a local coordinate system, such that, the orthogonal curve to the rotation is a planar extremal curve of the total curvature type energy $\mathbf{\Theta}_a$, \eqref{3}, with $a=\sqrt{-K}$. Therefore, all rotational surfaces verifying that $K<0$ is constant can be locally described as binormal evolution surfaces with planar critical curves of $\mathbf{\Theta}_{\sqrt{-K}}$ as initial condition and with velocity given by $\lvert \mathcal{I}\rvert$, see \eqref{i}.

Let $S$ be a rotational surface with constant negative curvature, $K<0$. Then, there exists a planar curve $\gamma$ in $\mathbb R^3$ such that $S$ can be seen as the set $\{\phi_t(\gamma(s))\, ; \, t\in\mathbb R\}$, where now $\{\phi_t\, ;\, t\in\mathbb R\}$ represents a one-parameter group of rotations. Therefore, from now on we are going to use the notation $S=S_\gamma$, and we are going to say that $\gamma$ is the \emph{profile curve}. Then, we have

\begin{thm}[\cite{P}] \label{thm2}
Let $S_\gamma$ be a rotational surface with profile curve $\gamma$ and such that its Gaussian curvature is constant and negative, that is, $K<0$ is constant. Then, the curvature of $\gamma$ verifies the Euler-Lagrange equation for $\mathbf{\Theta}_{\sqrt{-K}}$, \eqref{3}.
\end{thm}

\begin{figure}[H]
\begin{centering}
{\includegraphics[width=5cm,height=5cm]{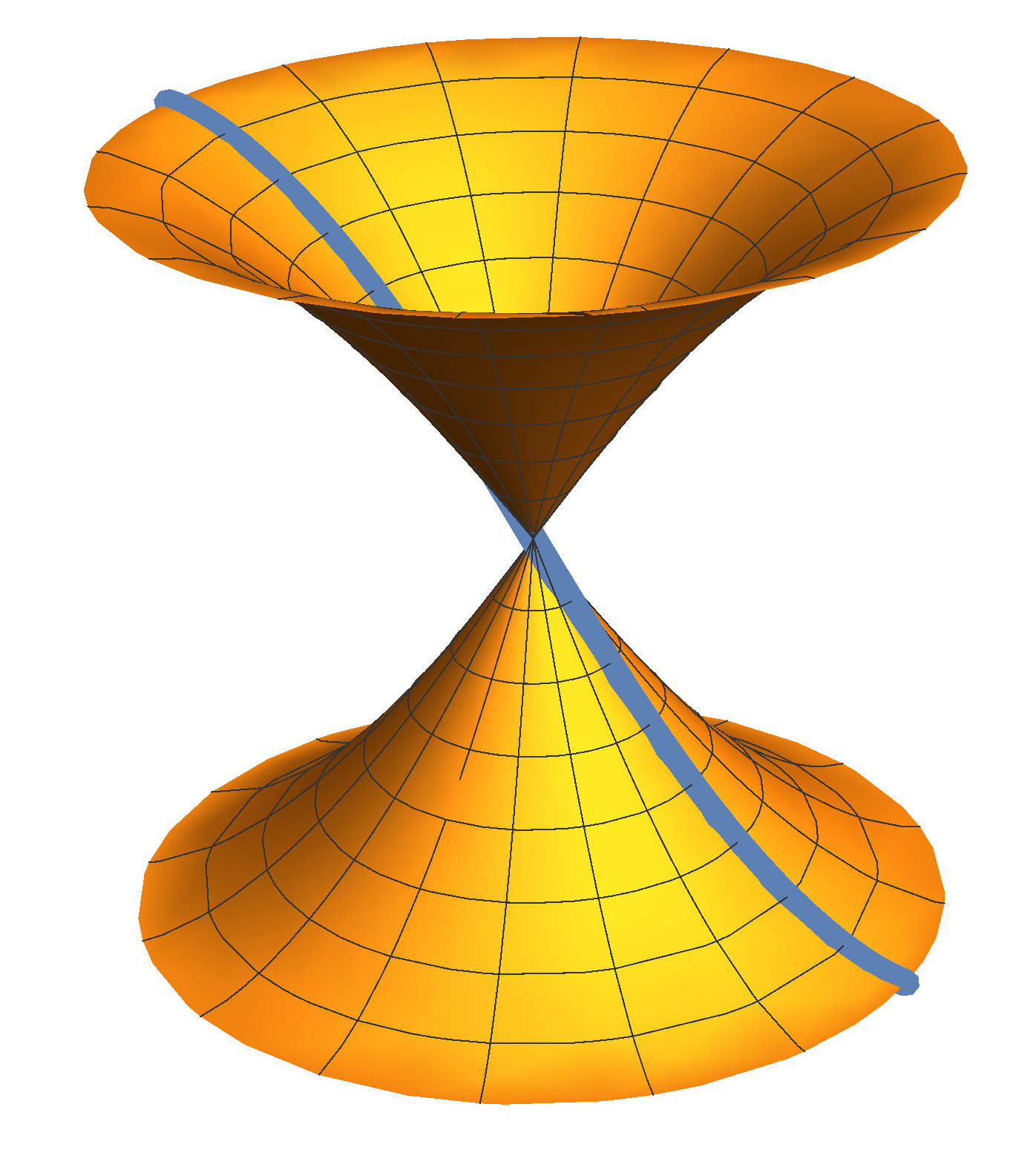}}\hspace{1cm}{\includegraphics[width=5cm,height=5cm]{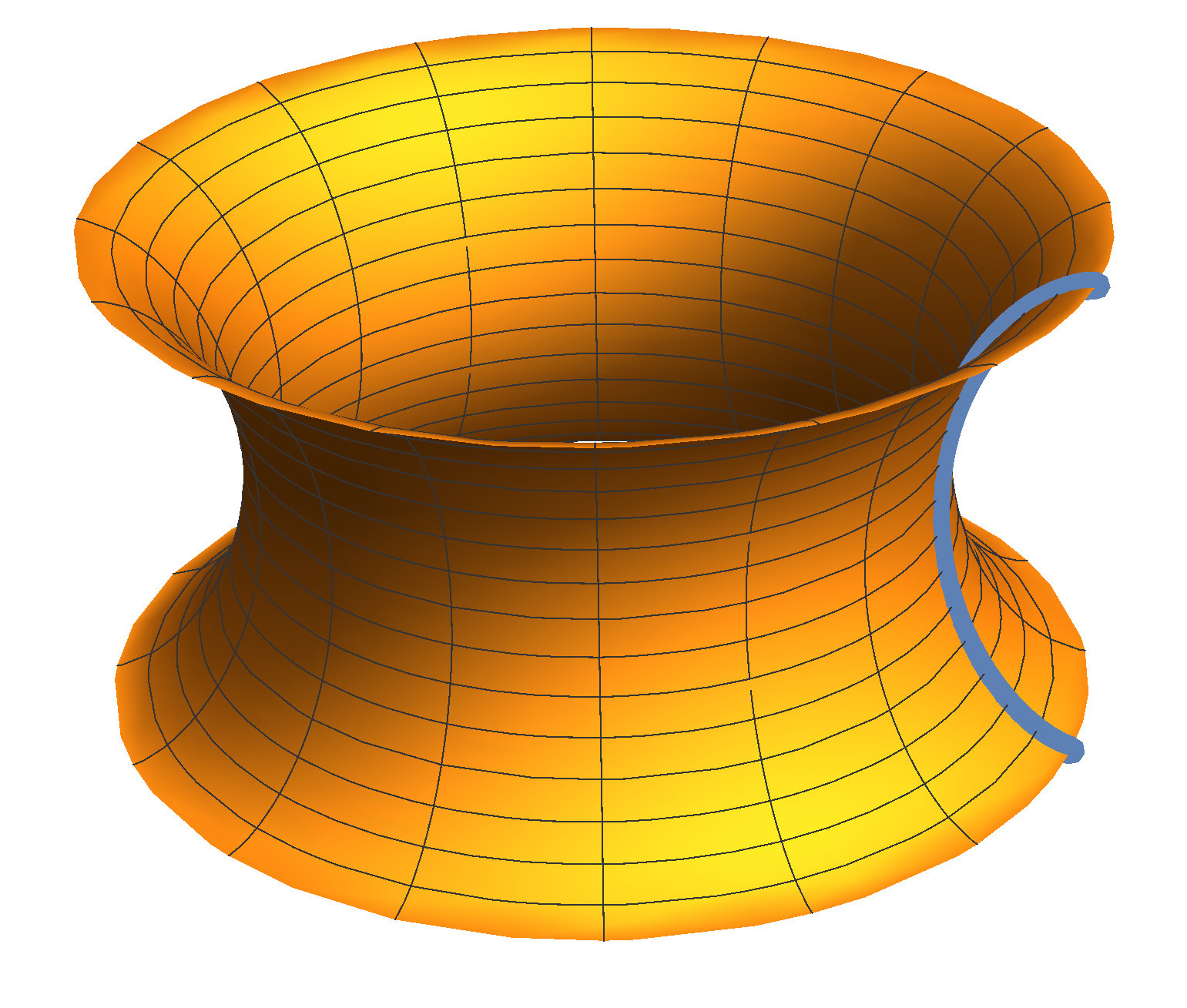}}\hspace{0.75cm}{\includegraphics[width=5cm,height=5cm]{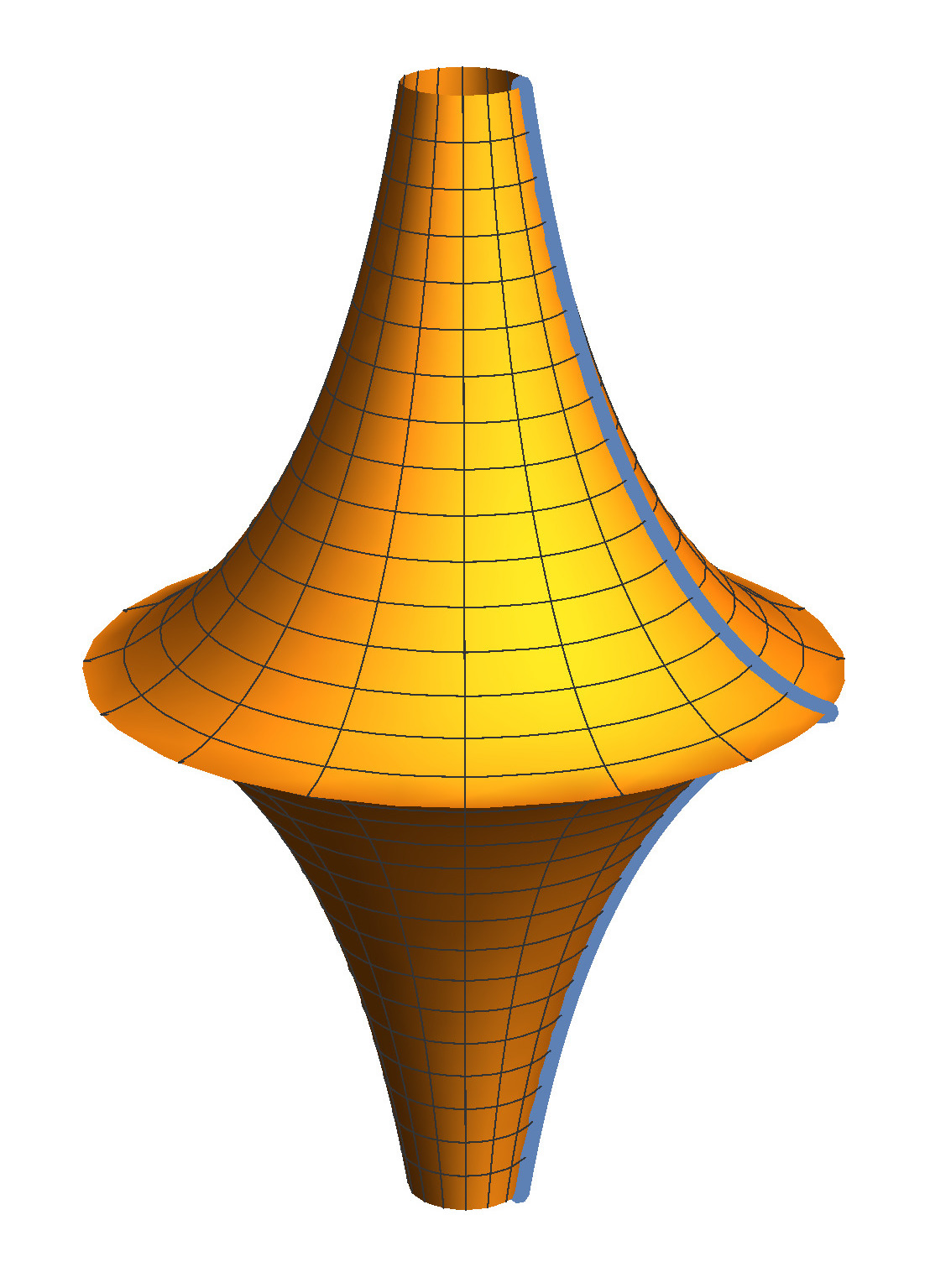}}
\caption{Rotational surfaces of $\mathbb{R}^3$ with constant negative Gaussian curvature and their profile curves. From left to right: conic type surface, hyperboloid type surface and pseudosphere.}\label{figure4}
\end{centering}
\end{figure}

Finally, for the sake of completeness, we summarize the complete description of all rotational surfaces with constant negative curvature and of their corresponding profile curves (see Figure \ref{figure4}). For a fixed positive real constant $a> 0$, consider that $K=-a^2$ is a negative constant, then the initial conditions of the binormal evolution have the curvatures given in Proposition \ref{curvaturas}, and they are defined respectively for the corresponding values of the arc-length parameter for any $d>a^2$;
\begin{enumerate}
\item If $f(z)=\sinh z$, then
\begin{equation}
s\in\left(\frac{-1}{a}\arcsinh\sqrt{\frac{a^2}{d-a^2}},\frac{1}{a}\arcsinh\sqrt{\frac{a^2}{d-a^2}}\right)\,,\nonumber
\end{equation}
and the binormal evolution surface is a rotational surface with constant negative Gaussian curvature of \emph{conic type}.
\item If $f(z)=\cosh z$, we also have that $d<2\,a^2$ and
\begin{equation}
s\in\left(\frac{-1}{a}\arccosh\sqrt{\frac{a^2}{d-a^2}},\frac{1}{a}\arccosh\sqrt{\frac{a^2}{d-a^2}}\right)\,.\nonumber
\end{equation}
Moreover, the binormal evolution surface swept out by this curve is a rotational constant negative Gaussian curvature surface of \emph{hyperbolic type}.
\item And, finally, if $f(z)=e^z$, then $s\in\left(-\infty,\frac{1}{a}\log\sqrt{\frac{a^2}{d-a^2}}\right)$, and the profile curve generates a rotational constant negative Gaussian curvature surface which is called \emph{pseudosphere}.
\end{enumerate}

Notice that none of these arc-length parameters is defined for all possible values of $s$. This implies that their corresponding binormal evolution surfaces are not complete, as expected by the statement of Hilbert's Theorem.

\section{Conclusions}

As mentioned in the introduction, there exist a one-to-one correspondence between sub-Riemannian geodesics of the unit tangent bundle of the plane, $\mathbb{R}^2\times\mathbb{S}^1$, used for visual curve completion and rotational surfaces of $\mathbb{R}^3$ with constant negative curvature.

In fact, as shown in \S 2, these sub-Riemannian geodesics can be characterized in terms of their projections into the plane, which turn out to be total curvature type extremals. Moreover, considering the plane, $\mathbb{R}^2$, as a totally geodesic surface of the Euclidean 3-space, $\mathbb{R}^3$, we have shown the existence of some kind of Noether's symmetries, which are described in the paper as Killing vector fields along the extremal curves. This can be regarded as a result coming from the nature of the variational problem.

In particular, one of these Killing vector fields along extremal curves is in the direction of the binormal to the curve. After considering its unique extension as a Killing vector field to the whole ambient space, it seems natural to consider the binormal evolution surface generated by evolving the extremal curve under its associated Killing vector field. Then, as shown in Theorem \ref{thm1}, we get a rotational surface of $\mathbb{R}^3$ with constant negative Gaussian curvature. On the other hand, the way back is also true as the statement of Theorem \ref{thm2} shows. Therefore, we have the desired one-to-one correspondence. 

This correspondence suggests that the sub-Riemannian structure of the unit tangent bundle of the plane used to model the first layer V1 may give some extra information. Indeed, first of all, we get the sub-Riemannian geodesic that completes the damaged or covered image, but, at the same time, we also get some rotational surfaces of $\mathbb{R}^3$ with constant negative curvature, which are closely connected to these geodesics.

Moreover, these rotational surfaces of $\mathbb{R}^3$ may be mapped to some special surfaces of $\mathbb{R}^2\times\mathbb{S}^1$, and, perhaps, the understanding of the properties of these related surfaces may give a new insight into the theory of visual curve completion.

\section*{Acknowledgements}

Research partially supported by MINECO-FEDER grant PGC2018-098409-B-100 and Gobierno Vasco grant IT1094-16. The author has been supported by Programa Posdoctoral del Gobierno Vasco, 2018. The author would like to thank Professor J. Arroyo for the images of Figure 1 and 2.

\end{document}